# A Smooth Curve as a Fractal Under the Third Definition

Ding Ma and Bin Jiang

Faculty of Engineering and Sustainable Development, Division of GIScience
University of Gävle, SE-801 76 Gävle, Sweden
Email: ding.ma|bin.jiang@hig.se

**Abstract**
It is commonly believed in the literature that smooth curves, such as circles, are not fractal, and only non-smooth curves, such as coastlines, are fractal. However, this paper demonstrates that a smooth curve can be fractal, under the new, relaxed, third definition of fractal – *a set or pattern is fractal if the scaling of far more small things than large ones recurs at least twice.* The scaling can be rephrased as a hierarchy, consisting of numerous smallest, a very few largest, and some in between the smallest and the largest. The logarithmic spiral, as a smooth curve, is apparently fractal because it bears the self-similar property, or the scaling of far more small squares than large ones recurs multiple times, or the scaling of far more small bends than large ones recurs multiple times. A half-circle or half-ellipse and the UK coastline (before or after smooth processing) are fractal, if the scaling of far more small bends than large ones recurs at least twice.

## 1. Introduction

Distinct from Euclidean geometry for dealing with 'cold' and 'dry' shapes, fractal geometry is a new geometric framework for dealing with irregular, rough, or non-smooth shapes that are ubiquitously seen in nature (Mandelbrot 1982). Despite of numerous applications of fractal geometry in geography (e.g., Batty and Longley 1994, Lam and Cola 2002, Chen 2008), its fundamental thinking has not been well adopted in cartography, with a few exceptions (Jiang et al. 2013, Jiang 2015b, Jiang 2017). For example, a cartographic curve is commonly understood as a set of more or similar line segments, which is the typical Euclidean geometric or non-recursive perspective (Jiang and Brandt 2016). Many cartographers are constrained by the non-recursive perspective or the traditional definitions of fractal, as recognized by Jiang and Yin (2014). A new definition of fractal has since then been offered as such: *a set or pattern is fractal if the scaling of far more small things than large ones recurs multiple times or with the ht-index being at least three* (Jiang and Yin 2014, Jiang 2015a, Gao et al. 2017). A cartographic curve should be more correctly perceived as a set of far more small bends than large ones, and the small bends are embedded in the large ones recursively (see Figure 2 for an illustration); a bend is determined by three vertices, but it is defined recursively, somehow like in Douglas algorithm (Douglas and Peucker 1973). This is the fractal geometric or recursive perspective.

Curves are often categorized as either smooth or non-smooth. Smooth curves are curves with a tangent at every point, or of functions with derivatives. On the contrary, non-smooth curves are curves without tangents, or of functions without derivatives. Non-smooth curves include those naturally occurring, such as coastlines and a mountain's profile or cross-section, and mathematically defined, such as the Koch curve and the Hilbert curve, which are examples of what is called the space-filling curves (Bader 2013). It is commonly accepted, in the literature, that a smooth curve, such as a circle, is unlikely to be fractal because it lacks self-similarity. As claimed by Mandelbrot (1982, p. 1), *"clouds are not spheres, mountains are not cones, coastlines are not circles, and bark is not smooth, nor does* *lightning travel in a straight line."* Fractal geometry provides a general framework for studying many irregular objects or non-smooth curves in nature (Falconer 2003). However, a smooth curve could be self-similar (Irving and Segerman 2013, Figure 2 therein in particular). The logarithmic spiral looks smooth, but is fractal because of its property of self-similarity. Under the new definition of fractal (cf. the next section), a set or pattern as a whole is iteratively divided into the head and the tail, with the head being recursively as a sub-whole. The self-similarity recurs between the head and the whole data or its sub-wholes, i.e., the head is recursively self-similar to the whole or sub-wholes. This paper demonstrates that a smooth curve can be a fractal because, if it is curved enough. We examined a series of curves including a half-circle, half-ellipse, the logarithmic spiral, and the UK coastline, to support our argument that a smooth curve can be a fractal under this new definition.

The remainder of this paper is structured as follows. Section 2 presents three definitions of fractal: the first two are based on power law relationship between the scale and the detail, while the third is relaxed from the power law, and based on head/tail breaks (Jiang 2013). Section 3 further illustrates head/tail breaks and the third definition using a simple curve with 10 recursively defined bends. Section 4 reports several experiments in order to show that smooth curve can be fractal under the third definition. Section 5 further discusses on the implications of this study for cartography, and finally Section 6 draws a conclusion.

**2. Three definitions of fractal**
Traditionally, fractals are defined from the top down, either strictly (definition 1) or statistically (definition 2) (Mandelbrot 1982, Cattani and Ciancio 2016). By the top down, we mean that a fractal is usually generated from a simple Euclidean shape, and by following a rule iteratively or endlessly; see below the Koch curve for example.

*Definition 1 and 2: A set or pattern is fractal if there is a power-law relationship between detail in the fractal (y) and the scale (x) at which it is measured, i.e., $y = x \wedge \alpha$, where α is the power-law exponent or the fractal dimension.*

A simple shape (such as a line segment or a square) is iteratively differentiated, being transformed into a complex shape with many levels of scale by following some simple generating rules. For example, the Koch curve is generated or transformed from a one-unit segment with the scale decreasing by one third: 1, 1/3, 1/9, 1/27, …, and with the detail increasing by four: 1, 4, 16, 64, …, leading to the power law relationship of $y = x \wedge -1.26$. The Koch curve is a good example of Definition 1 – strict or classic fractals, which are based on a strict power-law relationship, implying that all pairs of point (1, 1), (1/3, 4), (1/9, 16), (1/27, 64), …, are exactly on the distribution line $y = x \wedge -1.26$ (Figure 1a). The first definition is limited to classic fractals such as the Koch curve, the Hilbert curve, and Cantor set. The first definition was relaxed by Mandelbrot (1982), who defined so called statistical fractals – the second definition. It is relaxed in the sense that all pairs of point $(1 + \varepsilon_1, 1 + d_1)$, $(1/3 + \varepsilon_2, 4 + d_2)$, $(1/9 + \varepsilon_3, 16 + d_3)$, $(1/27 + \varepsilon_4, 64 + d_4)$, …, (where $\varepsilon_i$, and $d_i$ indicate small epsilons and deviations) are not exactly on but around the distribution line $y = x \wedge -1.26$ (Figure 1b). The second definition means that many real-world phenomena such as mountains, trees, clouds, and coastlines are fractal. This power-law relationship requirement, either strictly (definition 1) or statistically (definition 2), is too tough for many real-world fractals (Clauset et al. 2009), particularly, for fractals at earlier phases of development (Jiang and Yin 2014, Figure 2 therein). A new, relaxed definition of fractals does not require a power-law relationship between the scale and the detail. Instead, the so called third definition only requires recurring scaling pattern of far more small things than large ones (Jiang and Yin 2014, Jiang 2015a, Gao et al. 2017).

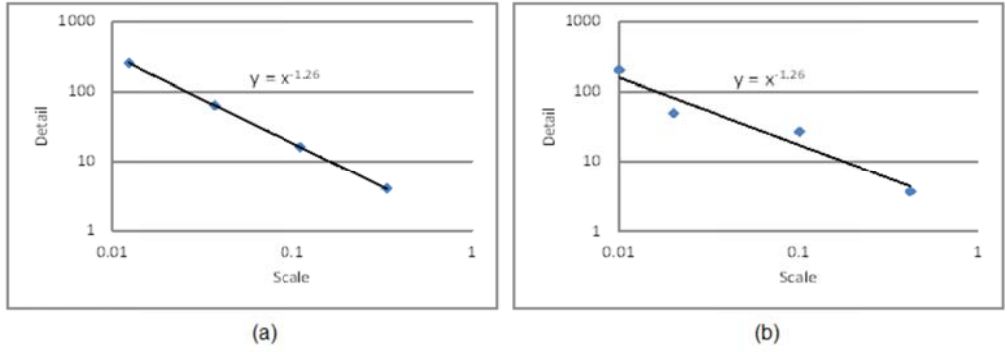

Figure 1: (Color online) Power-law relationship required for the first two definitions of fractal (Note: These two plots show a power law relationship, either strictly (a) or statistically (b), between the first four scales (1, 1/3, 1/9 and 1/27) or ($1 + \varepsilon_1$, $1/3 + \varepsilon_2$, $1/9 + \varepsilon_3$, and $1/27 + \varepsilon_4$) and the four details (1, 4, 16, 64) or ($1 + d_1$, $4 + d_2$, $16 + d_3$, $64 + d_4$) of the Koch curves. The difference between the two plots or between strictly and statistically is that with (a) the four points are exactly on the trend line, while with (b) the four points are around the trend line.)

*Definition 3: A set or pattern is fractal if the scaling of far more small things than large ones recurs multiple times or with the ht-index being at least three.*

This new and third definition (Jiang and Yin 2014) can help characterize fractals at different phases of development, with different degrees of complexity. For example, the Koch curve at iteration 3 is less complex than that at iteration 4. Unlike the first two definitions which are top-down, this new definition is from the bottom up. Given a whole set or pattern, it is iteratively divided into the head (being a sub-whole recursively) and the tail, in order to assess how many times the division of the head and the tail can occur, or to be more precise, how many times the head is self-similar to the whole or sub-whole. The number of times plus one, termed as the ht-index (Jiang and Yin 2014), indicates the complexity of the whole or the whole dataset.

It is important to note that the second definition is relaxed from the first, and the third definition is further relaxed from the second. Thus every new definition is inclusive to the previous one. Under the first definition, only the Koch curve is fractal, neither a coastline nor a highway is fractal. Under the second definition, both the Koch curve and a coastline are fractal, but a highway is not, since it is usually perceived to be smooth or less irregular than a coastline. However, under the third definition, not only the Koch curve and a coastline, but also a highway is fractal.

**3. Illustration of head/tail breaks and the third definition**
Let us examine these 10 numbers that strictly follow Zipf's law (1949): 1, 1/2, 1/3, …, and 1/10. The average of these 10 numbers is 0.29, which divides these numbers into two parts. The first three numbers, which are larger than the average, are called the head. The remaining seven, which are smaller than the average, are called the tail. The average for the head or the three largest numbers is 0.61, which divides the first three numbers into two parts again. The first number is larger than the second average, so is again the head. The remaining two are smaller than the second average, so are again the tail. This recursive division process is called the head/tail breaks (Jiang 2013), which leads to the ht-index that is defined by one plus the number of times that the scaling of far more small numbers than large ones recurs (Jiang and Yin 2014). The ht-index is actually the number of classes or hierarchical levels. Formally, the head/tail breaks is a recursive function:

```
Recursive Function Head/tail Breaks
   Break a data series, as a whole, into the head and the tail;
   // the head for those greater than the mean
   // the tail for those less than the mean
   While (head <= 40%): // 40% indicates a small head, and long tail
      Head/tail Breaks (head);
End Function
Ht-index = Number of iterations + 1
```

As we can see, the scaling of far more small numbers than large ones for the 10 numbers recurs twice, so the ht-index is three. Therefore, these 10 numbers, as a data set, are fractal. In order to further illustrate the third definition of fractal, we deliberately create a curve consisting of 10 bends (Figure 2), whose sizes ($x_i$) were exactly equal to the 10 numbers (Note: a bend is defined recursively, consisting of three vertices). Clearly, there are far more small bends than large ones, i.e., $x_1 + x_2 + x_3 > x_4 + x_5 + x_6 + x_7 + x_8 + x_9 + x_{10}$, and $x_1 > x_2 + x_3$. Therefore, the curve is considered to be a fractal.

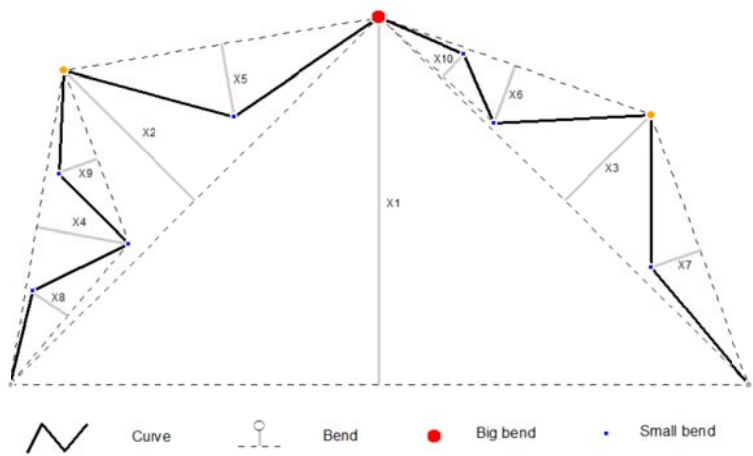

Figure 2: (Color online) Illustration of the new definition of fractal

(Note: The curve consists of 10 bends, which are imposed recursively, and scaling of far more small bends than large ones recurs twice: $x_1 + x_2 + x_3 > x_4 + x_5 + x_6 + x_7 + x_8 + x_9 + x_{10}$, and $x_1 > x_2 + x_3$, so the ht-index = 3. It should be noted that small bends are embedded in the large ones recursively. For example, $x_2$ and $x_3$ can be considered to be embedded in $x_1$.)

## 4. Experiments

We conducted several experiments on a set of smooth curves to illustrate the fact that a curve can be fractal if it is sufficiently curved, with the ht-index being at least three for its bends. The set of curves include a half-circle, two half-ellipses, the logarithmic spiral (Figure 3), and the UK coastline (Figure 4). We first partitioned these curves into numerous bends as illustrated in Figure 2, and then examine whether the scaling of far more small bends than large ones recurs multiple times. We found that these curves are all fractal under the new, relaxed, third definition. However, neither a half-circle nor half-ellipse is traditionally considered fractal. In what follows, we provide detailed evidence to support our argument that a smooth curve can be fractal.

### 4.1. Experiments: Case study I

We first created a half-circle with different numbers of vertices that are equally distributed along the half-circle curves: 128, 250, 500, 1000, and 6000. We then examined if the scaling of far more small bends than large ones recurs multiple times. For the first two cases, the ht-index is 4, which means that the scaling of far more small bends than large ones recurs three times. For the three other cases, the ht-index is 5, which is only slightly different from the first two cases. For all five cases, the bends follow a power-law distribution, with alpha approximately 1.61, $p \geq 0.65$. This power-law detection method is based on minimum likelihood, arguably the most rigorous power-law detection method (Clauset et al. 2009).

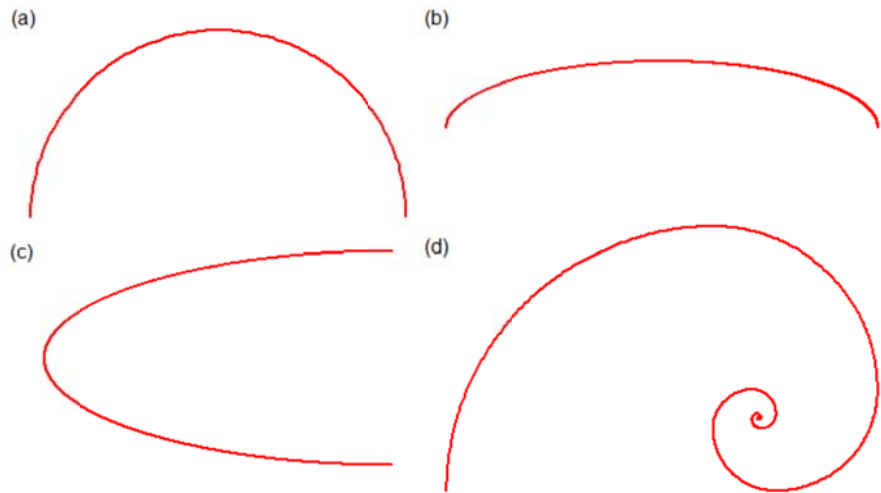

Figure 3: (Color online) Four smooth curves that are fractal under the third definition
(Note: (a) a half-circle, (b) a half-ellipse, (c) a half-ellipse, and (d) the logarithmic spiral)

**4.2. Experiments: Case study II**
From the first case study, we learned that the number of vertices has little effect on the ht-index. In the second study, we created an ellipse, and took its two halves with the same 5,998 vertices. The bends derived from the two half-ellipses have an ht-index of 4, implying that the scaling of far more small bends than large ones recurs three times. More importantly, the bends follow a striking power-law distribution, with alpha approximately 1.51, and $p \geq 0.01$.

**4.3. Experiments: Case study III**
We created five versions of the logarithmic spiral, also termed equiangular spiral or growth spiral, with different numbers of vertices: 37, 74, 138, 300, and 720, which are equally distributed along the curve. The first version did not have enough vertices for us to effectively verify it was fractal. The second and third versions of the logarithmic spiral have ht-indices of 4, implying that the scaling of far more small bends than large ones recurs three times. The fourth and fifth versions have ht-indices of 5, implying that the scaling of far more small bends than large ones recurs four times. In summary, the logarithmic spiral is fractal, with an ht-index of at least 4, given a reasonably sufficient number of vertices. The bends of the last four versions strikingly follow a power-law distribution, with an alpha ≈ 1.6, and $p \geq 0.12$. It is important to note that the logarithmic spiral often appears in nature, such as in a nautilus shell, a Romanesco broccoli, and a ram's horn (Thompson 1917).

The fact that the logarithmic spiral is fractal can be inferenced simply from its self-similarity property. However, both its fractal dimension and topological dimension are 1.0, which can also be seen in Table 1, so the logarithmic spiral is not fractal. Given the controversy, whether the logarithmic spiral is fractal has two opposite answers in the current literature (Chen 2017). In addition, the golden spiral is a special case of the logarithmic spiral, and the former can be approximately represented by the Fibonacci spiral, which comes from the Fibonacci sequence of 1, 1, 2, 3, 5, 8, 13, 21, 34, 55, 89, ..., whose scaling ratio of the subsequent number to the previous one always approaches to the golden ratio phi approximately 1.618. Statistically, there are far more small numbers than large ones in the Fibonacci sequence. Geometrically, there are far more small squares than large ones. All this evidence points to the fact that the logarithmic spiral is fractal.

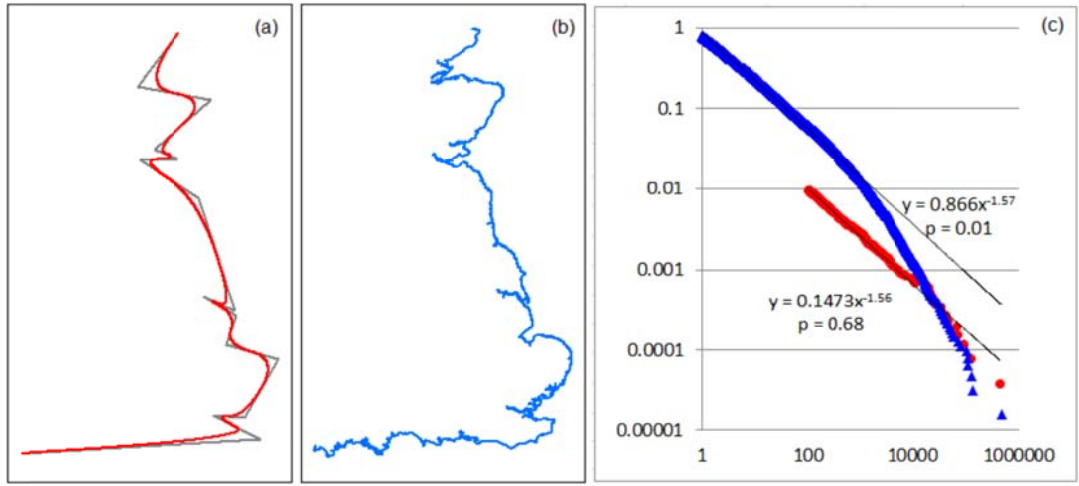

Figure 4: (Color online) The UK coastline before (gray) and after smoothing (red and blue)
(Note: The coastline with 14 vertices (gray), and 25,612 vertices (red) after smoothing (a); the coastline with 10,859 vertices (gray, but covered by the blue line), and with 62,550 vertices (blue) after smoothing (b); power-law distributions of the bends for the two versions of the smoothed UK coastline (c). It should be noted that due to a very high density of coastline vertices in Panel (b) (1) the smoothed coastline (blue) does not visually smooth, but is as smooth as the one in Panel (a) if one zoomed into it, and (2) the unsmoothed coastline (gray) is completely covered by the smoothed one.)

### 4.4. Experiments: Case study IV

We used the UK coastline for this case study, with 10,859 vertices as offered in the OpenStreetMap database (www.osm.org). We further smoothed the coastline using the popular Bezier interpolation method, so it becomes a curve with 62,550 vertices (Figure 4). For the same curve, before and after the smoothing operation, the ht-index is 7 (Table 1), implying that the scaling of far more small bends than large ones recurs six times. To further verify this, we used the same coastline with only 14 vertices (Figure 4, Table 1), for which the ht-index is 3, implying that the scaling of far more small bends than large ones recurs twice. After the smoothing operation, the curve of 14 vertices becomes a smooth curve with 25,612 vertices. This smoothed curve had an ht-index of 5, implying that the scaling of far more small bends than large ones recurs four times.

Table 1: Some of the curves are measured by bends and boxes
(Note: From the box-counting fractal dimension viewpoint, the logarithmic spiral is not fractal, but is fractal from the ht-index viewpoint. From the fractal dimension viewpoint, the coastline with 14 vertices before and after smoothing is fractal, but they remain fractal under the new definition or from the ht-index viewpoint.)

| Curves | Ht-index | Power law Alpha | p | Fractal dimension |
|---|---|---|---|---|
| Logarithmic spiral (300) | 5 | 1.58 | 0.25 | 0.987 |
| Coastline (10,859) | 7 | 2.15 | 0.83 | 1.203 |
| Smoothed coastline (62,550) | 7 | 1.56 | 0.68 | 1.202 |
| Coastline (14) | 3 | 2.52 | 0.02 | 1.115 |
| Smoothed coastline (25,612) | 5 | 1.57 | 0.01 | 1.097 |

The above experiments were conducted under the third definition of fractal. In addition, we used the box-counting technique for computing their fractal dimensions. The half-circle and the half-ellipses had a fractal dimension at range of [1.0, 1.1], which means they are not fractal under the second definition, since their fractal dimensions are close to 1.0. This is the same for the logarithmic spiral, which has a fractal dimension of 0.987 (Table 1). However, the fractal dimensions of the coastline with different numbers of vertices, either before or after smoothing, are very different from 1.0 (Table 1), implying that the smooth coastline is not a fake fractal – or a "fractal rabbit" (Kaye 1989, Chen 2015), but a true fractal, even under the traditional statistical definition – the second definition.

## 5. Implications

This study or the new definition of fractal in general implies that Euclidean geometric shapes, such as circles, are just a special case of fractals. Under the third definition of fractal, all shapes – regular or irregular, smooth or non-smooth – are fractal as long as the scaling of far more small things than large ones recurs multiple times. As we saw in the above experiments, a half-circle and a half-ellipse are fractal, as they are considered as a set of recursively defined bends: far more small bends than large ones, and small bends are embedded recursively in large ones. The small and large bends constitute respectively the tail and the head during the head/tail breaks process. The head is recursively self-similar to the whole or sub-wholes under the new definition of fractal. This study opens up new horizons for better understanding cartographic curves: smooth or non-smooth. They are regarded as a set of far more small bends than large ones – the fractal geometric perspective, rather than a set of more or less similar segments – the Euclidean geometric perspective. To further illustrate (Figure 5), the half circle is visualized by a set of bends of different colors: the blue bends with the smallest ht-indices, the red bend with the largest ht-index, and other color bends with ht-indices between the smallest and largest. Eventually, the generalization or different levels of scale of the half circle are no more than recursively keeping large bends or equivalently removing small ones. Note that at every level of scale, there are far more small bends than large ones.

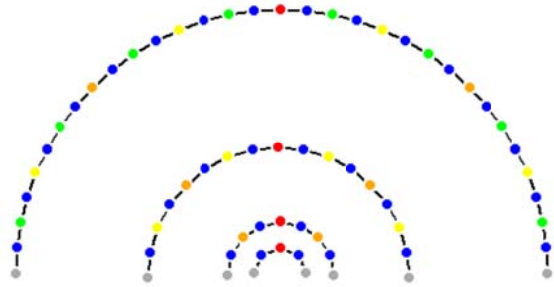

Figure 5: (Color online) A half circle at four different levels of scale
(Note: The spectral colors ranging from blue to red are used to indicate the ht-indices of individual bends – blue for the smallest, red for the largest, and other colors for those between the smallest and largest. At every level of scale, there are far more small bends than large ones, indicating fractal property of the half circle.)

## 6. Conclusion

Many smooth curves tend to be fractal under the new, relaxed, third definition of fractal, for the scaling of far more small bends than large ones recurs multiple times. This study points to the fact that a half-circle or a half-ellipse is fractal, as is a smoothed coastline. It should be made explicitly clear that a cartographic curve that is fractal under the second definition remains fractal after it has been smoothed, still under the second definition; see the column fractal dimension of Table 1, where fractal dimension remains very close before and after smoothing. The new technique for measuring hierarchical levels or the complexity in general – the number of times the scaling of far more small bends than large ones recurs is natural or organic. It is not like the mechanistic box-counting method for measuring fractal dimension. The new definition is from the bottom up, while the traditional classic or statistical definition is essentially top-down. Therefore, the new definition of fractal is more useful for better understanding real-world fractals or geographic features in particular.

## References:

Bader M. (2013), *Space-Filling Curves: An introduction with applications in scientific computing*, Springer: Berlin.

Batty M. and Longley P. (1994), *Fractal Cities: A geometry of form and function*, Academic Press: London.

Cattani C. and Ciancio A. (2016), On the fractal distribution of primes and prime-indexed primes by the binary image, *Physica A*, 460, 222-229.

Chen Y. G. (2008), *Fractal Urban Systems: Scaling, symmetry, and spatial complexity* (in Chinese), Science Press: Beijing.

Chen Y. G. (2015), Power-law distributions based on exponential distributions: Latent scaling, spurious Zipf's law, and fractal rabbits, *Fractals*, 23(2): 1550009.

Chen Y. G. (2017), Fractal analysis based on hierarchical scaling in complex systems, In: Fernando Brambila (ed), *Fractal Analysis - Applications in Health Sciences and Social Sciences*, Rijeka: InTech, 141–164.

Clauset A., Shalizi C. R., and Newman M. E. J. (2009), Power-law distributions in empirical data, *SIAM Review*, 51, 661-703.

Douglas D. and Peucker T. (1973), Algorithms for the reduction of the number of points required to represent a digitized line or its caricature, *The Canadian Cartographer*, 10(2), 112–122.

Falconer K. (2003), *Fractal geometry: Mathematical foundations and applications* (2nd edition), John Wiley & Sons Ltd: Chichester, England.

Gao P., Liu Z. Liu G., Zhao H., and Xie X. (2017), Unified Metrics for Characterizing the Fractal Nature of Geographic Features, *Annals of American Association of Geographers*, 1 – 17, http://dx.doi.org/10.1080/24694452.2017.1310022

Irving G. and Segerman H. (2013), Developing fractal curves, *Journal of Mathematics and the Arts*, 7(3-4), 103-121.

Jiang B. (2013), Head/tail breaks: A new classification scheme for data with a heavy-tailed distribution, *The Professional Geographer*, 65 (3), 482 – 494.

Jiang B. (2015a), Head/tail breaks for visualization of city structure and dynamics, *Cities*, 43, 69-77.

Jiang B. (2015b), The fractal nature of maps and mapping, *International Journal of Geographical Information Science*, 29(1), 159-174.

Jiang B. (2017), Line simplification, Richardson D., Castree N., Goodchild M. F., Kobayashi A., Liu W., and Marston R. A. (editors, 2017), *The International Encyclopedia of Geography*, John Wiley & Sons: New Jersey, 1-7, DOI: 10.1002/9781118786352.wbieg0005

Jiang B. and Brandt A. (2016), A fractal perspective on scale in geography, *ISPRS International Journal of Geo-Information*, 5(6), 95; doi:10.3390/ijgi5060095.

Jiang B. and Yin J. (2014), Ht-index for quantifying the fractal or scaling structure of geographic features, *Annals of the Association of American Geographers*, 104(3), 530–541.

Jiang B., Liu X. and Jia T. (2013), Scaling of geographic space as a universal rule for map generalization, *Annals of the Association of American Geographers*, 103(4), 844 – 855.

Kaye B. H. (1989), *A Random Walk through Fractal Dimensions*, VCH Publishers: New York.

Lam N. S-N, and Cola L. D. (2002), *Fractals in Geography*, The Blackburn Press: New Jersey.

Mandelbrot B. B. (1982), *The Fractal Geometry of Nature*, W. H. Freeman and Co.: New York.

Thompson D. W. (1917), *On Growth and Form*, Cambridge University Press: Cambridge.

Zipf G. K. (1949), *Human Behaviour and the Principles of Least Effort*, Addison Wesley: Cambridge, MA.